\documentclass[a4paper,12pt]{amsart}
\usepackage{amsmath,amsthm,amssymb}
\usepackage{mathtools}
\usepackage{ifthen}
\nonstopmode \numberwithin{equation}{section}
\newtheorem{thm}{Theorem}

\newtheorem{cor}{Corollary}
\newtheorem{lem}{Lemma}

\newtheorem{conj}{Conjecture}

\theoremstyle{definition}
\newtheorem{defn}{Definition}[section]

\newtheorem{prob}[equation]{Problem}

\newenvironment{rem}{%
\bigskip
\noindent \textsl{{\sl Remark. }}}{\bigskip}
\newenvironment{rems}{%
\bigskip
\noindent \textsl{{\sl Remarks. }}}{\bigskip}

\newcounter {own}
\def\theown {\thesection       .\arabic{own}}

\newenvironment{pf}[1][]{%
 \vskip 3mm
 \noindent
 \ifthenelse{\equal{#1}{}}%
  {{\slshape Proof. }}%
  {{\slshape #1.} }%
 }%
{\qed\bigskip}

\newcounter{alphabet}

\newcommand{\A}{{\mathcal A}}

\newcommand{\ID}{{\mathbb D}}

\newcommand{\IC}{{\mathbb C}}

\newcommand{\D}{{\mathbb D}}

\def\be{\begin{equation}}
\def\ee{\end{equation}}

\newcommand{\bee}{\begin{enumerate}}
\newcommand{\eee}{\end{enumerate}}

\newcommand{\blem}{\begin{lem}}
\newcommand{\elem}{\end{lem}}
\newcommand{\bthm}{\begin{thm}}
\newcommand{\ethm}{\end{thm}}
\newcommand{\bcor}{\begin{cor}}
\newcommand{\ecor}{\end{cor}}
\newcommand{\beg}{\begin{examp}}
\newcommand{\eeg}{\end{examp}}
\newcommand{\begs}{\begin{examples}}
\newcommand{\eegs}{\end{examples}}
\newcommand{\bdefe}{\begin{defn}}
\newcommand{\edefe}{\end{defn}}
\newcommand{\bprob}{\begin{prob}}
\newcommand{\eprob}{\end{prob}}
\newcommand{\bei}{\begin{itemize}}
\newcommand{\eei}{\end{itemize}}

\newcommand{\bcon}{\begin{conj}}
\newcommand{\econ}{\end{conj}}
\newcommand{\bcons}{\begin{conjs}}
\newcommand{\econs}{\end{conjs}}
\newcommand{\bprop}{\begin{propo}}
\newcommand{\eprop}{\end{propo}}
\newcommand{\br}{\begin{rem}}
\newcommand{\er}{\end{rem}}
\newcommand{\brs}{\begin{rems}}
\newcommand{\ers}{\end{rems}}
\newcommand{\bo}{\begin{obser}}
\newcommand{\eo}{\end{obser}}
\newcommand{\bos}{\begin{obsers}}
\newcommand{\eos}{\end{obsers}}
\newcommand{\bpf}{\begin{pf}}
\newcommand{\epf}{\end{pf}}
\newcommand{\ba}{\begin{array}}
\newcommand{\ea}{\end{array}}
\newcommand{\beq}{\begin{eqnarray}}
\newcommand{\beqq}{\begin{eqnarray*}}
\newcommand{\eeq}{\end{eqnarray}}
\newcommand{\eeqq}{\end{eqnarray*}}

\newcounter{minutes}
\divide\time by 60
\newcounter{hours}
\multiply\time by 60 \addtocounter{minutes}{-\time}
\begin{document}
\title{On the radius of concavity for certain classes of functions}
\begin{center}
{\tiny \texttt{FILE:~\jobname .tex,
        printed: \number\year-\number\month-\number\day,
        \thehours.\ifnum\theminutes<10{0}\fi\theminutes}
}
\end{center}

\author{Bappaditya Bhowmik${}^{~\mathbf{*}}$}
\address{Bappaditya Bhowmik, Department of Mathematics,
Indian Institute of Technology Kharagpur, Kharagpur - 721302, India.}
\email{bappaditya@maths.iitkgp.ac.in}
\author{Souvik Biswas}
\address{Souvik Biswas, Department of Mathematics, Indian Institute of Technology Kharagpur, Kharagpur - 721302, India.}
\email{souvikbiswas158@gmail.com}


\begin{abstract}
  Let $\mathcal{A}$ denote the class of all analytic functions $f$ defined in the open unit disc $\D$ with the normalization $f(0)=0=f'(0)-1$ and let $P'$ be the class of functions $f\in\mathcal{A}$ such that ${\rm{Re}}\,f'(z)>0$, $z\in\D$. In this article, we obtain radii of concavity of $P'$ and for the class $P'$ with the fixed second coefficient. After that, we consider linearly invariant family of functions, along with the class of starlike functions of order $1/2$ and investigate their radii of concavity. Next, we obtain a lower bound of radius of concavity for the class of functions $\mathcal{U}_0(\lambda)=~\{f\in\mathcal{U}(\lambda) : f''(0)=0\}$, where
  $$
  \mathcal{U}(\lambda)=\left\{f\in\mathcal{A} : \left|\left(\frac{z}{f(z)}\right)^2f'(z)-1\right|<\lambda,~z\in \D\right\},\quad \lambda \in (0,1].
  $$
  We also investigate the meromorphic analogue of the class $\mathcal{U}(\lambda)$ and compute its radius of concavity. 
\end{abstract}
\thanks{The first author of this article would like to thank
NBHM, India for its financial support through a research grant (Ref. No.- 02011/10/2022/R\&D-II/14948).}

\maketitle
\pagestyle{myheadings}
\markboth{B. Bhowmik and S. Biswas}{On the radius of concavity for certain classes of functions}

\maketitle

\bigskip

\section{Introduction and Preliminaries}
  Throughout this article, we will use the following notations. Let $\mathbb{C}$ be the whole complex plane and $\bar{\mathbb{C}}:=\mathbb{C}\cup\{\infty\}$ be the Riemann sphere. We denote the open unit disc of the complex plane by $\D$, i.e. $\D:=\{z\in\mathbb{C} : |z|<1\}$. Let $\mathcal{A}$ be the class of all analytic functions in $\D$, satisfying the normalization $f(0)=0=f'(0)-1$.  Let $\mathcal{S}$ be the class of all univalent functions in $\A$. For years, researchers have been exploring different subclasses of $\mathcal{S}$ that possess geometric properties. Among these, the most prominent classes include the class of convex functions of order $\beta$, $0\leq\beta<1$, defined as 
 $$
  C(\beta)=\{f\in \mathcal{A} : {\rm{Re}}\,\left(1+\frac{zf''(z)}{f'(z)}\right)>\beta, ~ z\in\D\}
 $$
 and the class of all starlike functions of order $\beta$, $0\leq\beta<1$, defined as
 $$
  S^\ast(\beta)=\{f\in \mathcal{A} : {\rm{Re}}\,\left(\frac{zf'(z)}{f(z)}\right)>\beta, ~ z\in \D\}.
 $$
 We denote the class of all convex functions by $C:=C(0)$, which consists of all functions $f\in\mathcal{S}$ such that $f$ maps $\D$ conformally onto a convex domain. Likewsie, the class of starlike functions is denoted by $S^\ast:=S^\ast(0)$, which consists of all functions $f \in \mathcal{S}$ such that $f$ maps $\D$ conformally onto a starlike domain with respect to the origin. We first recall here that, the {\it{radius of convexity}} (or {\it{starlikeness}}) for a subset $\mathcal{A}_1$ of $\mathcal{A}$, is the largest value of $r \in (0,1]$ such that every function $f \in \mathcal{A}_1$ is convex (or starlike) within the disk $\mathbb{D}_r$. In $1920$, Nevanlinna (c.f. \cite{nev}) proved that the radius of convexity for $\mathcal{S}$ is $2-\sqrt{3}$. Later, in $1934$, Grunsky obtained the radius of starlikeness for $\mathcal{S}$ as $\tanh{\pi/4}$ (see \cite[p.~141]{grunsky}). In this article, our main concern is to find {\it{radius of concavity}} for certain classes of analytic and meromorphic univalent functions. More precisely, in \cite{sou}, we introduced the notion of {\it{radius of concavity}} in the context of analytic and meromorphic mappings separately. We briefly recall these definitions from \cite{sou}. Let $Co(A)$ be the class of functions $f$ that satisfy the following conditions:
 \begin{itemize}
     \item[(i)] $f\in\mathcal{S}$ with the additional condition $f(1)=\infty$.
     \item[(ii)] $\mathbb{C}\setminus f(\D)$ is convex.
     \item[(iii)] The opening angle of $f(\D)$ at $\infty$ is less than or equal to $\pi A$, $A\in (1,2]$.
 \end{itemize}
 In \cite[Theorem~2]{fg}, Avkhadiev and Wirths proved that $f \in Co(A)$ if and only if
 $$
 {\rm{Re}}\,T_f(z)>0, \quad \forall ~z \in \D,
 $$
 where, $f(0)=0=f'(0)-1$ and
\begin{equation}\label{P2eq1.1}
T_f(z):=\frac{2}{A-1}\left[\frac{(A+1)}{2}\frac{1+z}{1-z}-1-z\frac{f''(z)}{f'(z)}\right].
\end{equation}
In the analytic case, the {\it{radius of concavity}} (with respect to $Co(A)$) for a subclass of $\A$ is defined as below (see \cite[Definition~$1.2$]{sou}).
\begin{defn}\label{P2defn1.1}
The radius of concavity (with respect to $Co(A)$) of a subset $\A_1$ of $\A$ is the largest number $R_{Co(A)}\in (0,1]$ such that for each function $f \in \A_1$, ${\rm{Re}}\,T_f(z)>0$ for all $|z|<R_{Co(A)}$, where $T_f$ is defined in \eqref{P2eq1.1}.
\end{defn}
For meromorphic functions, the definition of radius of concavity is not the same as above. Again we recall it from \cite{sou}. Let $\mathcal{A}(p)$ be the class of meromorphic functions in $\D$ with a simple pole at $z=p$, $p\in(0,1)$ and normalized by the condition $f(0)=0=f'(0)-1$. Let $S(p)$ denote the set of all univalent functions in $\mathcal{A}(p)$. Let $Co(p)$ be the class of functions $f\in S(p)$ such that $\overline{\IC}\setminus f(\ID)$ is a bounded
convex set. In $1971$, J. A. Pfaltzgraff and B. Pinchuk (see \cite [p.~145]{pfaltz}) proved that $f\in Co(p)$ if and only if $f\in S(p)$ such that
$$
{\rm{Re}}\,P_f(z)>0, \quad  z\in \ID, \quad P_f(p)=\frac{1+p^2}{1-p^2} \quad {\rm{and}} \quad P_f(0)=1,
$$
where
\begin{equation}\label{P2eq1.2}
P_f(z):=-\left[1+\frac{zf''(z)}{f'(z)}+\frac{z+p}{z-p}-\frac{1+pz}{1-pz}\right].
\end{equation}
For meromorphic functions, the {\it{radius of concavity}} (with respect to $Co(p)$) for a subclass of $\A(p)$ is defined as below (see \cite[Definition~$1.1$]{sou}).
\begin{defn}\label{P2defn1.2}
The radius of concavity (with respect to $Co(p)$) of a subset $\A_1(p)$ of $\A(p)$  is the largest number $R_{Co(p)}\in(0,1]$
such that for each function $f \in \A_1(p)$, ${\rm{Re}}\,P_f(z)>0$ for all $|z|<R_{Co(p)}$, where, $P_f$ is defined in \eqref{P2eq1.2}.
\end{defn}
We mention here that, in \cite[Theorem~$4$]{sou} we obtained the radius of concavity of $\mathcal{S}$. We also derived a lower bound for the radius of concavity of $S(p)$ (see \cite[Theorem~$2$]{sou}). Furthermore, we determined the same for the linear combinations of functions belonging to $S(p)$, $\mathcal{S}$, and $Co(A)$ with complex coefficients.

Let $P'$ be the class of functions $f\in\A$ such that ${\rm{Re}}\,f'(z)>0$, $z\in\D$. In $1962$, Macgregor (c.f. \cite{mac}) determined the radius of convexity of the class $P'$ as $\sqrt2-1$ with an extremal function $f_0(z)=-z-2\ln{(1-z)}$. In $2004$, Todorov (c.f. \cite{to}) gave an alternative proof of this result and generalized all the extremal functions associated with it. Todorov also considered the class $P'$ with fixed second coefficient and obtained the radius of convexity for the same. In this article, we will compute the radius of concavity of $P'$ using Definition~\ref{P2defn1.1}. Furthermore, we will determine a lower bound for the radius of concavity of the class $P'$ with fixed second coefficient.

Let $\mathcal{LS}$ be the set of locally univalent functions defined in $\D$, with the normalization $f(0)=0=f'(0)-1$. Let ${\rm{Aut}}(\D)$ denote the set of holomorphic automorphisms of $\D$. If $f \in \mathcal{LS}$ and $\phi\in{\rm{Aut}}(\D)$, then the Koebe transform of $f$ with respect to $\phi$ is defined as 
$$
 \Lambda_\phi(f)(z):=\frac{(f\circ\phi)(z)-(f\circ\phi)(0)}{(f\circ\phi)'(0)}, \quad z\in\D.
$$
We recall that a family $\mathcal{F}$ is called a {\it{linear-invariant family}} (L.I.F.) if $\mathcal{F} \subset \mathcal{LS}$ and $\Lambda_\phi(f) \in \mathcal{F}$ for all $f \in \mathcal{F}$, $\phi \in {\rm{Aut}}(\D)$. The order of a L.I.F. $\mathcal{F}$ is 
$$
 {\rm{ord}}~\mathcal{F}:=\sup\left\{\left|\frac{f''(0)}{2}\right| : f \in \mathcal{F}\right\}.
$$
For a detailed study of this linear-invariant family, we urge the reader to go through \cite[Ch.~$5$]{graham}. In $1964$, Pommerenke (c.f. \cite{Po}) proved that the radius of convexity of a linear-invariant family $\mathcal{F}$ is $\alpha-\sqrt{\alpha^2-1}$, where ${\rm{ord}}~\mathcal{F}=\alpha<\infty$. In this article, we will compute the radius of concavity of a linear-invariant family $\mathcal{F}$. As every convex function is starlike, a function that is convex in $\D$, satisfies ${\rm{Re}}\,\left(zf'(z)/f(z)\right)>0$ for $z\in\D$. In fact, for $f\in C$, the aforementioned condition can be improved to ${\rm{Re}}\,\left(zf'(z)/f(z)\right)>1/2$ for $z\in\D$ (see \cite{str}). In other words, if $f \in C$, then $f \in S^\ast(1/2)$. In \cite{Mac}, Macgregor obtained the radius of convexity of the class $S^\ast(1/2)$ as $\sqrt{2\sqrt{3}-3}$. In this article, we will determine a lower bound for the radius of concavity of $S^\ast(1/2)$.
 
In the final section of this article, we investigate radius of concavity for another interesting subclass $\mathcal{U}(\lambda)$, $\lambda \in(0,1]$ of $\mathcal{S}$. We first present a brief overview about this class of functions. For $\lambda\in(0,1]$, let
$$
\mathcal{U}(\lambda)=\{f \in \mathcal{A} : |U_f(z)|<\lambda, ~ z\in\D\},
$$
where $U_f(z):=(z/f(z))^2f'(z)-1$. We know that $\mathcal{U}(\lambda) \subset \mathcal{S}$, for $0<\lambda \leq 1$, (c.f. \cite{oz}). For more details about this class, we urge the reader to go through the articles \cite{ob, oz, po} and references therein. Let $\mathcal{U}_0(\lambda)$ be the class of functions $f$ in $\mathcal{U}(\lambda)$ with the additional condition $f''(0)=0$. In \cite{va}, Ponnusamy and Vasundhra proved that $\mathcal{U}_0(\lambda) \subset S^\ast$ for $0<\lambda \leq 1/\sqrt2$ and $\mathcal{U}_0(\lambda) \subset S^\ast(1/2)$ for $0<\lambda \leq 1/3$. They also proposed a conjecture that each $f\in \mathcal{U}_0(\lambda)$ is convex in $\D$ if $0<\lambda \leq 3-2\sqrt2$. In \cite{vas}, Vasundhra also derived a lower bound for the radius of convexity of $\mathcal{U}(\lambda)$ and $\mathcal{U}_0(\lambda)$. In \cite{kar}, it is shown that the above conjecture is not valid and an improved lower bound for the radius of convexity of $\mathcal{U}_0(\lambda)$ has been found. In this article, we will determine a lower bound for the radius of concavity of the class $\mathcal{U}_0(\lambda)$. In  \cite{par}, the first author of this article and F. Parveen  considered a meromorphic analog of the class $\mathcal{U}(\lambda)$ as follows. For $\lambda\in(0,1]$, let
$$
 \mathcal{V}_p(\lambda)=\{f \in \A(p) : |U_f(z)|<\lambda,~ z\in\D\}.
$$
We urge the reader to go through the articles \cite{par2, par} for more information about this class of functions. Finally, we will obtain a lower bound for the radius of concavity of the class $\mathcal{V}_p(\lambda)$. As a byproduct of this result, we also find radius of convexity for functions in the class $\mathcal{V}_p(\lambda)$. It is worth to mention here that, in certain cases, the obtained radii may not be sharp; and determining the exact radii will be interesting problems of study in future.

\section{Radius of concavity of functions whose derivative has a positive real part}
We now state and prove the first result of this article.
\begin{thm}
    The radius of concavity $R_{Co(A)}$ of $P'$ is 
    \begin{equation}\label{P2eq2.1}
         1-\frac{2}{\sqrt{A+3}}, \quad A\in(1,2].
    \end{equation}
This radius is sharp.
\end{thm}
\begin{pf}
       Since $f\in P'$, therefore $f'\prec(1+z)/(1-z)$ and this implies
        \begin{equation}\label{P2eq2.2}
              f'(z)=\frac{1+g(z)}{1-g(z)}, \quad z\in\D,
        \end{equation}
        where $g$ is an analytic function in $\mathbb{D}$ with $g(0)=0$ and $|g(z)|<1$. From \eqref{P2eq2.2}, we get
        \begin{equation}\label{P2eq2.3}
          1+\frac{zf''(z)}{f'(z)}=1+\frac{2zg'(z)}{1-(g(z))^2}, \quad z\in\D.
        \end{equation}
        By the Schwarz-Pick lemma, we get
        $$
            |g'(z)| \leq \frac{1-|g(z)|^2}{1-|z|^2}, \quad z\in\D. 
        $$
        Using \eqref{P2eq2.3} and the above inequality, we get
        \begin{equation}\label{P2eq2.4}
         {\rm{Re}}\,\left(1+\frac{zf''(z)}{f'(z)}\right) \leq 1+\frac{2r}{1-r^2}, \quad |z|=r<1.
        \end{equation}
        It is easy to see that
        \begin{equation}\label{P2eq2.5}
         {\rm{Re}}\,\left(\frac{1+z}{1-z}\right)\geq \frac{1-r}{1+r}, \quad |z|=r<1.
        \end{equation}
        Applying \eqref{P2eq2.4} and \eqref{P2eq2.5}, we get
        $$
         {\rm{Re}}\,T_f(z)\geq\frac{2}{A-1}\left(\left(\frac{A+1}{2}\right)\frac{1-r}{1+r}-1-\frac{2r}{1-r^2}\right), \quad |z|=r<1.
        $$
        The right hand side quantity of the above inequality is strictly positive if
        $$
         r^2-2r+\frac{A-1}{A+3}>0.
        $$
        Hence,
        ${\rm{Re}}\,T_f(z)>0$ if $|z|=r<R_{Co(A)}$,
        where, $R_{Co(A)}$ is given in \eqref{P2eq2.1}. Also, if we consider $f_0(z)=-z+2\log{(1+z)}$, $z\in \D$, then for this function we compute 
        $$
         T_{f_0}(z)= \left(\frac{A+3}{A-1}\right)\frac{1}{1-z^2}\left(z^2+2z+\frac{A-1}{A+3}\right).
        $$
        We observe that, if $z=-r$ and $R_{Co(A)}<|z|<1$, then ${\rm{Re}}\,T_{f_0}(z)<0$. This proves the sharpness of the radius of concavity.
 \end{pf}
 
Let $P'(a)$ be a subclass of $P'$, consisting of analytic functions with the expansion
$$
 f(z)=z+az^2+\sum_{n=3}^\infty a_nz^n, \quad z\in\D,
$$
having the fixed second coefficient $a=f''(0)/2$. It is known that $0\leq a \leq 1$ (see \cite[Section~4]{to}). In the following theorem, we determine a lower bound for the radius of concavity of $P'(a)$.
 \begin{thm}
     If $f \in P'(a)$, then ${\rm{Re}}\,T_f(z)>0$ for $|z|<R_{Co(A)}$, where, $R_{Co(A)}$ is the least value of $r\in (0,1)$ satisfying $\phi(r)=0$ with 
     \begin{equation*}
      \begin{split}
        \phi(r):= & r^4-\frac{2(A+1)(1-a)}{A+3}r^3+\frac{2(A(1-2a)-2a-3)}{A+3}r^2\\
                  & -\frac{2(A(1-a)+3a+1)}{A+3}r+\frac{A-1}{A+3}.
      \end{split}
     \end{equation*}
 \end{thm}
 \begin{pf}
      From \cite{to}, we have 
 $$
     \left|\frac{zf''(z)}{f'(z)}\right| \leq \frac{2|z|}{1-|z|^2}\left(\frac{a|z|^2+2|z|+a}{|z|^2+2a|z|+1}\right), \quad z\in\D, \quad 0 \leq a \leq 1.
 $$
 From the above inequality, we get
 $$
  {\rm{Re}}\, \left(1+\frac{zf''(z)}{f'(z)}\right) \leq 1+\frac{2|z|}{1-|z|^2}\left(\frac{a|z|^2+2|z|+a}{|z|^2+2a|z|+1}\right),\quad z\in\D.
 $$
 Applying \eqref{P2eq2.5} and the above inequality, we get
 $$
  {\rm{Re}}\,T_f(z) \geq \frac{2}{A-1}\left(\left(\frac{A+1}{2}\right)\frac{1-r}{1+r}-1-\frac{2r}{1-r^2}\left(\frac{ar^2+2r+a}{r^2+2ar+1}\right)\right), \quad |z|=r<1.
 $$
 The right hand side of the above inequality is strictly positive if $|z|<R_{Co(A)}$, where, $R_{Co(A)}$ is given in the statement of the theorem. We now investigate the existence of $R_{Co(A)}$ for each $A \in (1,2]$ and $a\in[0,1]$. We see that the function $\phi$ which is defined in the statement of the theorem, is continuous on $[0,1]$ with
$$
\phi(0)=\frac{A-1}{A+3}>0 \quad {\rm{and}} \quad \phi(1)=-8\left(\frac{a+1}{A+3}\right)<0.
$$
Therefore, by the intermediate value theorem, $\phi$ has at least one root in $(0,1)$. Hence, $R_{Co(A)}$ exists for every $A \in (1,2]$ and $a\in[0,1]$.
 \end{pf}

\section{Radius problems for linear-invariant family and for starlike functions of order $1/2$}
We start this section with certain results related to linear invariance and radius problems. In the following theorem, we determine the radius of concavity for a linear-invariant family $\mathcal{F}$ of order $\alpha<\infty$.
\begin{thm}
    The radius of concavity $R_{Co(A)}$ of a linear-invariant family $\mathcal{F}$ of order $\alpha<\infty$ is 
    \begin{equation}\label{P2eq3.1}
        \frac{A+1+2\alpha-2\sqrt{(A+\alpha)(1+\alpha)}}{A-1}.
    \end{equation}
    This radius is sharp.
\end{thm}
\begin{pf}
    For a linear-invariant family $\mathcal{F}$ of order $\alpha$, Pommerenke (c.f. \cite{Po}) proved that
    $$
     \alpha=\sup_{f\in \mathcal{F}} \sup_{|z|<1} \left|-\overline{z}+\frac{1}{2}(1-|z|^2)\frac{f''(z)}{f'(z)}\right|. 
    $$
    This implies that for each $f\in\mathcal{F}$,
    $$
     \left|-\overline{z}+\frac{1}{2}(1-|z|^2)\frac{f''(z)}{f'(z)}\right| \leq \alpha, \quad z\in \D.
    $$
    From the above inequality, we get
    $$
     \left|\frac{zf''(z)}{f'(z)}-\frac{2|z|^2}{1-|z|^2}\right|\leq \frac{2\alpha |z|}{1-|z|^2},\quad z\in\D.
    $$
    Thus, we get for each $f\in\mathcal{F}$,
    $$
     \left|\frac{zf''(z)}{f'(z)}-\frac{2r^2}{1-r^2}\right|\leq \frac{2\alpha r}{1-r^2}, \quad |z|=r<1.
    $$
    Therefore, we get
    $$
     {\rm{Re}}\,\left(\frac{zf''(z)}{f'(z)}\right) \leq \frac{2\alpha r+2r^2}{1-r^2}, \quad |z|=r<1.
    $$
    From \eqref{P2eq2.5} and the above inequality, we get
    $$
     {\rm{Re}}\,T_f(z) \geq \frac{2}{A-1}\left(\left(\frac{A+1}{2}\right)\frac{1-r}{1+r}-1-\frac{2\alpha r+2r^2}{1-r^2}\right), \quad |z|=r<1.
    $$
    The right hand side of the above inequality is strictly positive if
    $$
    r^2-2\left(\frac{A+1+2\alpha}{A-1}\right)r+1>0.
    $$
    Hence, 
    $$
     {\rm{Re}}\,T_f(z)>0 \quad {\rm{if}} ~ |z|=r<R_{Co(A)},
    $$
    where, $R_{Co(A)}$ is given in \eqref{P2eq3.1}. Also, if we consider 
    $$
     g_0(z)=\frac{1}{2\alpha}\left(1-\left(\frac{1-z}{1+z}\right)^{\alpha}\right), \quad z \in \D,
    $$ 
    then for this function, we compute 
        $$
         T_{g_0}(z)=\frac{1}{1-z^2}\left(z^2+2\left(\frac{A+1+2\alpha}{A-1}\right)z+1\right).
        $$
        We observe that, if $z=-r$ and $R_{Co(A)}<|z|<1$, then ${\rm{Re}}\,T_{g_0}(z)<0$. This proves sharpness of the radius of concavity.
    
\end{pf}

\begin{rem}
    Since ${\rm{ord}}~\mathcal{S}=2$, therefore from \eqref{P2eq3.1}, we get
    $$
    R_{Co(A)}=\frac{A+5-\sqrt{12(A+2)}}{A-1},
    $$
    which is the radius of concavity for the class $\mathcal{S}$ (see \cite[Theorem~$4$]{sou}). 
\end{rem}

In the next theorem, we obtain the radius of concavity of $S^{\ast}(1/2)$. We need the following lemma by Macgregor (see \cite[Lemma~$1$]{Mac}) to prove our next theorem.
\begin{lem}\label{P2lem1}
    The function $g(z)$ is analytic for $|z|<1$ and satisfies $g(0)=1$ and ${\rm{Re}}\,g(z)>1/2$ for $|z|<1$ if and only if $g(z)=1/(1+z\phi(z))$, where $\phi(z)$ is analytic and satisfies $|\phi(z)|\leq 1$ for $|z|<1$.
\end{lem}
\begin{thm}
    If $f \in S^\ast (1/2)$, then ${\rm{Re}}\,T_f(z)>0$ for $|z|<R_{Co(A)}$, where, $R_{Co(A)}$ is the least value of $r \in (0,1)$ satisfying $u(r)=0$ with
    $$
     u(r)=(A-1)r^3-(3A+1)r^2+(3A+7)r-A-1.
    $$
\end{thm}
\begin{pf}
    If we consider the function $zf'/f$, then from the Lemma~\ref{P2lem1}, we have 
    \begin{equation}\label{P2eq3.2}
          \frac{zf'(z)}{f(z)}=\frac{1}{1+z\phi(z)},
    \end{equation}
    where $\phi$ is analytic in $\D$ satisfying $|\phi(z)| \leq 1$ for $z\in \D$. From the Schwarz-Pick lemma, we get 
    \begin{equation}\label{P2eq3.3}
     |\phi'(z)| \leq \frac{1-|\phi(z)|^2}{1-|z|^2}, \quad z \in \D.
    \end{equation}
    Differentiating \eqref{P2eq3.2} and by simple calculations, we have 
    $$
     \frac{zf''(z)}{f'(z)} +1=\frac{1-z\phi(z)-z^2\phi(z)}{1+z\phi(z)}, \quad z \in \D.
    $$
    This implies
    $$
     \frac{zf''(z)}{f'(z)} +1=\frac{(1-z\phi(z)-z^2\phi(z))(\overline{1+z\phi(z)})}{|1+z\phi(z)|^2}, \quad z \in \D.
    $$
    Therefore,
    \begin{equation}\label{P2eq3.4}
         {\rm{Re}}\,\left(\frac{zf''(z)}{f'(z)} +1\right)=\frac{1-|z|^2|\phi(z)|^2}{|1+z\phi(z)|^2}-\frac{{\rm{Re}}\,\left(z^2\phi'(z)(1+\overline{z\phi(z)})\right)}{|1+z\phi(z)|^2}, \quad z\in\D.
    \end{equation}
    Since
    \beqq
    |1+z\phi(z)|^2 & = & 1+|z|^2|\phi(z)|^2+2 ~{\rm{Re}}\,(z\phi(z)) \\
                   & \geq & 1+|z|^2|\phi(z)|^2-2 ~|z||\phi(z)|=(1-|z||\phi(z)|)^2,
    \eeqq
    therefore, from \eqref{P2eq3.4} we get
    $$
        {\rm{Re}}\,\left(\frac{zf''(z)}{f'(z)} +1\right)\leq \frac{1-|z|^2|\phi(z)|^2}{(1-|z||\phi(z)|)^2}+\frac{\left|z^2\phi'(z)(1+\overline{z\phi(z)})\right|}{(1-|z||\phi(z)|)^2}, \quad z\in \D.
    $$
    With the help of \eqref{P2eq3.3}, the above inequality yields
    $$
         {\rm{Re}}\,\left(\frac{zf''(z)}{f'(z)} +1\right)\leq \frac{1-r^2|\phi(z)|^2}{(1-r|\phi(z)|)^2}+\frac{r^2(1-|\phi(z)|^2)(1+r|\phi(z)|)}{(1-r^2)(1-r|\phi(z)|)^2}
    $$
     for $|z|=r<1$.
    Applying \eqref{P2eq2.5} and the above inequality, we get 
    $$
     {\rm{Re}}\,T_f(z) \geq \frac{2}{A-1}\left(\frac{A+1}{2}~\frac{1-r}{1+r}-\frac{1-r^2|\phi(z)|^2}{(1-r|\phi(z)|)^2}-\frac{r^2(1-|\phi(z)|^2)(1+r|\phi(z)|)}{(1-r^2)(1-r|\phi(z)|)^2}\right)
    $$
    for $|z|=r<1$.
    If we let $|\phi(z)|=x$, $0\leq x\leq1$, then upon simplification, the above inequality yields 
    \begin{equation}\label{P2eq3.5}
    {\rm{Re}}\,T_f(z)\geq \frac{p(x)}{(A-1)(1-r^2)(1-xr)^2},\quad |z|=r<1,
    \end{equation}
   where,
    \begin{equation*}
       \begin{split}
         p(x)= & 2r^3 x^3+\left((A-1)r^2-2(A+1)r+A+5\right)r^2 x^2-2\big((A+2)r^2-2(A+1)r \\
               & +A+1\big)rx+(A+1)r^2-2(A+1)r+A-1.
       \end{split}
    \end{equation*}
    For $x\in[0,1]$, we get
     \begin{equation*}
        \begin{split}
             p'(x)=& 6r^3x^2+2r^2\left((A-1)r^2-2(A+1)+A+5\right)x \\
                         & -2r\left((A+2)r^2-2(A+1)r+A+1\right).
        \end{split}
    \end{equation*}
   It is a routine exercise to see that $p'$ is increasing in $[0,1]$ and therefore has the maximum value $p'(1)=2r\left((A-1)r^3-(3A+1)r^2+(3A+7)r-A-1\right)$. This implies that for $|z|=r<1$,
    $$
     p'(x) \leq 2ru(r), \quad 0 \leq x\leq 1,
    $$
     where, $u(r)=(A-1)r^3-(3A+1)r^2+(3A+7)r-A-1$, $r \in [0,1]$. As $u$ is continuous on $[0,1]$ with 
     $$
      u(0)=-(A+1)<0 \quad {\rm{and}} \quad u(1)=4>0,
     $$
     therefore, by the intermediate value theorem, $u$ has at least one root in $(0,1)$. Also $u$ is increasing in $[0,1]$ and has the minimum value $u(0)<0$. Let $R_{Co(A)}$ be the least value of $r \in (0,1)$ satisfying the equation $u(r)=0$. If $r<R_{Co(A)}$, then $u(r)<0$. This implies that $p'(x) < 0$ for $0\leq x \leq1$, if $r<R_{Co(A)}$. Therefore, if $0<r<R_{Co(A)}$ then $p$ is decreasing in $[0,1]$ and has the minimum value
    $$
      p(1)=(A+3)r^4-4(A+1)r^3+(6A+2)r^2-4(A+1)r+A+3=Q(r)~({\rm{say}}).
    $$ 
    Thus, we get from \eqref{P2eq3.5}
    $$
        {\rm{Re}}\,T_f(z)\geq \frac{Q(r)}{(A-1)(1-r^2)(1-xr)^2}, \quad |z|=r<R_{Co(A)}.
    $$
    By simple computations, it can be shown that $Q$ is decreasing on $[0,1]$ and has the minimum value $Q(1)=0$. This implies that $Q(r)>0$ for $r\in [0,R_{Co(A)})$. Therefore, ${\rm{Re}}\,T_f(z)>0$ if $|z|=r<R_{Co(A)}$, which we need to show.  
\end{pf}

\section{Radius problems for the classes $\mathcal{U}_0(\lambda)$ and $\mathcal{V}_p(\lambda)$}
In \cite{vas}, Vasundhra established a lower bound for the radius of convexity of $\mathcal{U}_0(\lambda)$. Subsequently, this lower bound has been improved in \cite[Theorem~2.4]{kar}. In the following theorem, we determine a lower bound for the radius of concavity of $\mathcal{U}_0(\lambda)$.  
\begin{thm}
   If $f \in \mathcal{U}_0(\lambda)$, then ${\rm{Re}}\ T_f(z)>0$ for $|z|<R_{Co(A)}$, where, $R_{Co(A)}=\min\{r_1,r_2\}$ with
   $$
   r_1=\sqrt{\frac{5+\lambda-\sqrt{(1-\lambda)(25-\lambda)}}{6\lambda}}
   $$
   and $r_2$ being the least positive value of $r$ satisfying $\phi(r)=0$ with
   $$
       \phi(r)=-\lambda(9-A)r^3-\lambda(A+11)r^2-(A+3)r+A-1.
   $$
\end{thm}
\begin{pf}
    For each $f \in \mathcal{U}_0(\lambda)$, we have from \cite[Theorem~$2.4$]{kar}
    $$
        \left|\frac{zf''(z)}{f'(z)}\right| \leq \frac{6\lambda r^2}{1-\lambda r^2},\quad |z|=r<r_1,
    $$
    where $r_1$ is given in the statement of the theorem. This implies that
 $$
     {\rm{Re}}\, \left(\frac{zf''(z)}{f'(z)}\right) \leq \frac{6\lambda r^2}{1-\lambda r^2}, \quad |z|=r<r_1.
 $$
 Applying \eqref{P2eq2.5} and the above inequality, we get
 $$
     {\rm{Re}}\, T_f(z) \geq \frac{2}{A-1}\left(\left(\frac{A+1}{2}\right)\frac{1-r}{1+r}-1-\frac{6\lambda r^2}{1-\lambda r^2} \right), \quad |z|=r<r_1.
 $$
 The right hand side of the above inequality is strictly positive if $|z|<r_2$, where $r_2$ is defined in the statement of the theorem. Thus, ${\rm{Re}}\,T_f(z)>0$ if $|z|<R_{Co(A)}$, where $R_{Co(A)}=\min\{r_1,r_2\}$. We now investigate the existence of $r_2$ for each $\lambda \in (0,1]$ and $A\in(1,2]$. We see that the function $\phi$ which is defined in the statement of the theorem, is continuous on $[0,1]$ with
$$
\phi(0)=A-1>0 \quad {\rm{and}} \quad \phi(1)=-20\lambda-4<0.
$$
Therefore, by the intermediate value theorem, $\phi$ has at least one root in $(0,1)$. Hence, this proves the existence of $r_2$  for every $\lambda \in (0,1]$ and $A\in(1,2]$.
\end{pf}

In the next theorem, we determine a lower bound for the radius of concavity of the class $\mathcal{V}_p(\lambda)$.
\begin{thm}
    If $f \in \mathcal{V}_p(\lambda)$, then ${\rm{Re}}\, P_f(z)>0$ for $|z|<R_{Co(p)}$, where, $R_{Co(p)}=\min\{r_1,r_2\}$ with
    $$
    r_1=\sqrt{\frac{3-\lambda-\sqrt{(1-\lambda)(9-\lambda)}}{2\lambda}}
    $$
    and $r_2$ being the least value of $r \in (0,p)$ satisfying $\phi(r)=0$ with 
    
\begin{equation*}
 \begin{split}
   \phi(r):= & (\lambda^2p^2)r^5+\lambda p(3\lambda p^2-\lambda+1)r^4+\lambda (5\lambda p^2-4p^2+3)r^3 \\
             & {+ p(1-4\lambda -\lambda p^2)r^2-(1+p^2+3\lambda p^2)r+p.}
 \end{split}
\end{equation*}

\end{thm}
\begin{pf}
We know that from \cite[Section~$2$]{par} if $f\in \mathcal{V}_p(\lambda)$, then there exists a holomorphic function $w_1$ such that $w_1(\D)\subseteq\overline{\D}$ and
\begin{equation}\label{P2eq4.1}
 \frac{z}{f(z)}=1-\left(\frac{f''(0)}{2}\right)z+\lambda z \int_0^z w_1(t) dt,\quad z\in \D.
\end{equation}
Let 
\begin{equation}\label{P2eq4.2}
    w(z):=\left(\int_p^zw_1(t)dt\right)/(z-p), \quad z\in \mathbb{D}.
\end{equation}
Then from \cite[Theorem~$1$]{par2}, we get 
\begin{equation}\label{P2eq4.3}
    \frac{z}{f(z)}=\frac{-(z-p)(1-\lambda pzw(z))}{p}, \quad z\in\D,
\end{equation}
where $w$ is analytic in $\mathbb{D}$ and $|w(z)|\leq 1$, $z \in \mathbb{D}$. By differentiating \eqref{P2eq4.1} with respect to $z$, we get
$$
 -z\left(\frac{z}{f(z)}\right)'+\frac{z}{f(z)}=1-\lambda z^2w_1(z), \quad {\rm{for}} ~ z\in\D.
$$
This implies that
$$
\left(\frac{z}{f(z)}\right)^2f'(z)=1-\lambda z^2w_1(z), \quad z\in\D.
$$
Then from \eqref{P2eq4.3}, we get
$$
    f'(z)=\frac{p^2(1-\lambda z^2 w_1(z))}{(z-p)^2(1-\lambda pzw(z))^2}, \quad z\in\D.
$$
Thus, we have for $z\in\D$,
\begin{equation}\label{P2eq4.4}
 \frac{zf''(z)}{f'(z)}=-\frac{2z}{z-p}-\frac{\lambda z^2(2w_1(z)+zw_1'(z))}{1-\lambda z^2w_1(z)}+\frac{2\lambda pz(w(z)+zw'(z))}{1-\lambda pzw(z)}, \quad z\in\D.
\end{equation}
From \eqref{P2eq4.2}, we get
$$
 (z-p)w'(z)+w(z)=w_1(z).
$$
Then by a simple calculation we see that 
$$
 zw'(z)+w(z)=\frac{1}{z-p}(zw_1(z)-pw(z)).
$$
By using this, we deduce from \eqref{P2eq4.4} that for $z\in\D$,
\begin{equation}\label{P2eq4.5}
    1+\frac{zf''(z)}{f'(z)}+\frac{z+p}{z-p}=-\frac{\lambda z^2(2w_1(z)+zw_1'(z))}{1-\lambda z^2w_1(z)}+\frac{2\lambda pz(zw_1(z)-pw(z))}{(z-p)(1-\lambda pzw(z))}.
\end{equation}
Let us define 
$$
 u(z)=z^2w_1(z), \quad z\in\mathbb{D}. 
$$
Then we get $u(0)=0=u'(0)$ and $|u(z)|\leq|z|^2$, $z \in \mathbb{D}$. The equation \eqref{P2eq4.5} becomes 
\begin{equation}\label{P2eq4.6}
    1+\frac{zf''(z)}{f'(z)}+\frac{z+p}{z-p}=-\frac{\lambda zu'(z)}{1-\lambda zu(z)}+\frac{2\lambda p(u(z)-pzw(z))}{(z-p)(1-\lambda pzw(z))}.
\end{equation}
Since $w_2(z)=(u(z))/z$ is an analytic function from $\D$ to $\D$, we get
$$
 |w_2'(z)|\leq \frac{1-|w_2(z)|^2}{1-|z|^2}, \quad z\in \D,
$$
which is equivalent to 
$$
 |zu'(z)-u(z)|\leq \frac{|z|^2-|u(z)|^2}{1-|z|^2}, \quad z\in \D.
$$
From \eqref{P2eq4.6} and the above inequality, we get
$$
 \left|1+\frac{zf''(z)}{f'(z)}+\frac{z+p}{z-p}\right| \leq \frac{\lambda (|z|^2-|u(z)|^2)}{(1-|z|^2)(1-\lambda |u(z)|)}+\frac{\lambda |u(z)|}{(1-\lambda |u(z)|)}+\frac{2\lambda p(|u(z)|+p|z|)}{||z|-p|(1-\lambda p |z|)},
$$
for $z\in\D$. If we let $|u(z)|=x$ and $|z|=r$ (note that $0 \leq x \leq r^2$), the above inequality becomes
\begin{equation}\label{P2eq4.7}
     \left|1+\frac{zf''(z)}{f'(z)}+\frac{z+p}{z-p}\right| \leq \frac{\lambda (r^2-x^2)}{(1-r^2)(1-\lambda x)}+\frac{\lambda x}{1-\lambda x}+\frac{2\lambda p(x+pr)}{(p-r)(1-\lambda p r)}.
\end{equation}
Let us denote 
$$
 \phi(x)= \frac{\lambda (r^2-x^2)}{(1-r^2)(1-\lambda x)}+\frac{\lambda x}{1-\lambda x}+\frac{2\lambda p(x+pr)}{(p-r)(1-\lambda p r)}, \quad {\rm{for}} \quad 0\leq x \leq r^2.
$$
Then by some simple computations, we see that
$$
 \phi(x) \leq \phi(r^2)=\frac{2\lambda r^2}{1-\lambda r^2}+\frac{2\lambda pr(p+r)}{(p-r)(1-\lambda pr)}, \quad {\rm{if}} \quad r<r_1,
$$
where $r_1$ is given in the statement of the theorem. Thus from \eqref{P2eq4.7}, we get
\begin{equation}\label{P2eq4.8}
    \left|1+\frac{zf''(z)}{f'(z)}+\frac{z+p}{z-p}\right| \leq \frac{2\lambda r^2}{1-\lambda r^2}+\frac{2\lambda pr(p+r)}{(p-r)(1-\lambda pr)}, \quad |z|=r<r_1.
\end{equation}
This implies that
\begin{equation}\label{P2eq4.9}
    {\rm{Re}}\, \left(1+\frac{zf''(z)}{f'(z)}+\frac{z+p}{z-p}\right) \leq \frac{2\lambda r^2}{1-\lambda r^2}+\frac{2\lambda pr(p+r)}{(p-r)(1-\lambda pr)}, \quad |z|=r<r_1.
\end{equation}
It is easy to see that
$$
 {\rm{Re}}\, \left(\frac{1+pz}{1-pz}\right) \geq \frac{1-pr}{1+pr}, \quad |z|=r<1.
$$
Applying \eqref{P2eq4.9} and the above inequality, we get
$$
 {\rm{Re}}\, P_f(z) \geq \frac{1-pr}{1+pr}-\frac{2\lambda r^2}{1-\lambda r^2}-\frac{2\lambda pr(p+r)}{(p-r)(1-\lambda pr)}, \quad |z|=r<r_1.
$$
The right hand side of the above inequality is strictly positive if $|z|=r<r_2$, where $r_2$ is defined in the statement of the theorem. Thus, ${\rm{Re}}\, P_f(z)>0$ if $|z|<R_{Co(p)}$, where, $R_{Co(p)}={\rm{min}}~\{r_1,r_2\}$. We now investigate the existence of $r_2$ for each $\lambda \in (0,1]$ and $p\in(0,1)$. The function $\phi$ which is defined in the statement of the theorem, is continuous on $[0,p]$ with
$$
\phi(0)=p>0 \quad {\rm{and}} \quad \phi(p)=-4\lambda p^3(1+p^2)(1-\lambda p^2)<0.
$$
Therefore, by the intermediate value theorem, $\phi$ has at least one root in $(0,p)$. This establishes the existence of $r_2$  for every $\lambda \in (0,1]$ and $p\in(0,1)$.
\end{pf}

In the following Corollary, we find radius of convexity for the class $\mathcal{V}_p(\lambda)$ by applying the estimate derived in $\eqref{P2eq4.8}$.
\begin{cor}
    If $f \in \mathcal{V}_p(\lambda)$, then $f$ maps $|z|<R_c$ onto a convex set, where, $R_c={\rm{min}}\{r_1,r_2\}$ with
    $$
     r_1=\sqrt{\frac{3-\lambda-\sqrt{(1-\lambda)(9-\lambda)}}{2\lambda}}
    $$
    and $r_2$ being the least value of $r \in (0,p)$ satisfying $\psi(r)=0$ with
    $$
     \psi(r)=(\lambda^2p)r^5+\lambda (1+2\lambda p^2)r^4-\lambda p(1-5\lambda p^2)r^3+(1-5\lambda p^2)r^2-p(2+3\lambda p^2)r+p^2.
    $$
\end{cor}
\begin{pf}
   From \eqref{P2eq4.8}, we have 
   $$
       {\rm{Re}}\, \left(1+\frac{zf''(z)}{f'(z)}\right) \geq {\rm{Re}}\,\left(\frac{p+z}{p-z}\right) -\frac{2\lambda r^2}{1-\lambda r^2}-\frac{2\lambda pr(p+r)}{(p-r)(1-\lambda pr)}, \quad |z|=r<r_1.
   $$
   It is easy to see that
   $$
    {\rm{Re}}\, \left(\frac{p+z}{p-z}\right) \geq \frac{p-r}{p+r}, \quad |z|=r<1.
   $$
Thus, we get 
$$
 {\rm{Re}}\, \left(1+\frac{zf''(z)}{f'(z)}\right) \geq \frac{p-r}{p+r} -\frac{2\lambda r^2}{1-\lambda r^2}-\frac{2\lambda pr(p+r)}{(p-r)(1-\lambda pr)}, \quad |z|=r<r_1.
$$
The right hand side of the above inequality is strictly positive if $|z|<r_2$, where $r_2$ is defined in the statement of the theorem. Thus, ${\rm{Re}}\, \left(1+zf''(z)/f'(z)\right)>0$ if $|z|<R_c$, where, $R_c={\rm{min}}\{r_1,r_2\}$. We now investigate the existence of $r_2$ for each $\lambda \in (0,1]$ and $p\in(0,1)$. The function $\psi$ which is defined in the statement of the theorem, is continuous on $[0,p]$ with
$$
\psi(0)=p^2>0 \quad {\rm{and}} \quad \psi(p)=-8\lambda p^4(1-\lambda p^2)<0.
$$
Therefore, by the intermediate value theorem, $\psi$ has at least one root in $(0,p)$. Thus, the existence of $r_2$ is confirmed for every $\lambda \in (0,1]$ and $p\in(0,1)$.
\end{pf}

\noindent{\bf Declarations of interest}: The authors have no relevant financial or non-financial interests to disclose.
The authors have no competing interests to declare that are relevant to the content of this article.\\
\noindent{\bf Data availability statement}: Data sharing not applicable to this article as no datasets were generated or analysed during the current study.


\begin{thebibliography}{99}

\bibitem{fg}{\sc F. G. Avkhadiev and K.-J. Wirths}: Concave schlicht functions with bounded opening angle at infinity, \textit{Lobachevskii
J. Math.} {\bf 17} (2005), 3--10.

\bibitem{sou}{\sc B. Bhowmik and S. Biswas}: Distortion, radius of concavity and several other radii results for certain classes of functions, \textit{Comput. Methods Funct. Theory}, (2024), DOI : https://doi.org/10.1007/s40315-024-00525-8.

\bibitem{par2}{\sc B. Bhowmik and F. Parveen}: On the Taylor coefficients of a subclass of meromorphic univalent functions, \textit{Bull. Malays. Math. Sci. Soc.} {\bf 42} (2019), no. 2, 793--802.

\bibitem{par}{\sc B. Bhowmik and F. Parveen}: Sufficient conditions for univalence and study of a class of meromorphic univalent functions, \textit{Bull. Korean Math. Soc.} {\bf 55} (2018), no. 3, 999--1006.

\bibitem{graham}{\sc I. Graham and G. Kohr}: Geometric function theory in one and higher dimensions, \textit{Monographs and Textbooks in Pure and Applied Mathematics} {\bf 255.} Marcel Dekker, Inc., New York (2003), {\rm{xviii}}+530.

\bibitem{grunsky}{\sc H. Grunsky}: Zwei Bemerkungen zur konformen Abbildung, \textit{Jber. Deutsch. Math.-Verein.} {\bf 43} (1934), 140--143. 

\bibitem{kar}{\sc V. Karunakaran and K. Bhuvaneswari}: On the radius of convexity for a class of conformal maps, \textit{Colloq. Math.} {\bf 109} (2007), no. 2, 251--256.

\bibitem{mac}{\sc Thomas H. Macgregor}: Functions whose derivative has a positive real part, \textit{trans. Amer. Math. Soc.} {\bf 104} (1962), 532--537.

\bibitem{Mac}{\sc Thomas H. Macgregor}: The radius of convexity for starlike functions of order ${1/2}$, \textit{Proc. Amer. Math. Soc.} {\bf 14} (1963), 71--76.

\bibitem{nev}{\sc R. Nevanlinna}: {\"U}ber die konforme Abbildung von Sterngebieten, \textit{\"Oversikt av Finska Vetenskaps-Soc. F\"orh.} {\bf 63(A)} (1920-21), 1--21.

\bibitem{ob}{\sc M. Obradovi\'c and S. Ponnusamy and K.-J. Wirths}: Geometric studies on the class $\mathcal{U}(\lambda)$, \textit{Bull. Malays. Math. Sci. Soc.} {\bf 39} (2016), no. 3, 1259--1284.

\bibitem{oz}{\sc S. Ozaki and M. Nunokawa}: The Schwarzian derivative and univalent functions, \textit{Proc. Amer. Math. Soc.} {\bf 33} (1972), 392--394.

\bibitem{pfaltz}{\sc J. Pfaltzgraff and B. Pinchuk}: A variational method for classes of meromorphic functions, \textit{J. Analyse Math.} {\bf 24} (1971), 101--150.

\bibitem{Po}{\sc C. Pommerenke}: Linear-invariante Familien analytischer Funktionen, {\rm{I}}, \textit{Math. Ann.} {\bf 155} (1964), no. 2, 108--154.

\bibitem{va}{\sc S. Ponnusamy and P. Vasundhra}: Criteria for univalence, starlikeness and convexity, \textit{Ann. Polon. Math.} {\bf 85} (2005), no. 2, 121--133.

\bibitem{po}{\sc S. Ponnusamy, K.-J. Wirths}: Coefficient problems on the class $\mathcal{U}(\lambda)$, \textit{Porbl. Anal. Issues Anal.} {\bf 7}(25) (2018), no. 1, 87--103.

\bibitem{str}{\sc E. Strohh$\ddot{a}$cker}: Beitr$\ddot{a}$ge zur Theorie der schlichten Funktionen, \textit{Math. Z.} {\bf 37} (1933), 356--380.

\bibitem{to}{\sc Pavel G. Todorov}: On the radius of convexity of functions whose derivative has a positive real part, \textit{J. Dyn. Syst. Geom. Theor.} {\bf 2} (2004), 35--41.

\bibitem{vas}{\sc P. Vasundhra}: Application of Hadamard product and hypergeometric functions in univalent function theory, \textit{Thesis}, Department of Mathematics, Indian Institute of Technology, Madras, June 2004.
    
\end{thebibliography}
\end{document}